\documentstyle[11pt]{article}

\newcommand{\qed}{\hfill$\Box$}

\newcommand{\inv}{{\it inv}}
\newcommand{\D}{{\it Des}}
\newcommand{\des}{{\it des}}
\newcommand{\M}{\stackrel{\longleftarrow}{\mbox{Min}}}
\newcommand{\m}{\stackrel{\longleftarrow}{\mbox{min}}}
\newcommand{\Der}{\stackrel{\longleftarrow}{\mbox{Min}}}
\newcommand{\der}{\stackrel{\longleftarrow}{\mbox{min}}}

\newcommand{\maj}{{\it maj}}

\newcommand{\fm}{{\it fmaj}}

\newcommand{\al}{\alpha}

\newtheorem{thm}{Theorem}[section]
\newtheorem{pro}[thm]{Proposition}
\newtheorem{lem}[thm]{Lemma}

\newtheorem{cor}[thm]{Corollary}
\newtheorem{fac}[thm]{Fact}

\newtheorem{exa}[thm]{Example}
\newtheorem{df}[thm]{Definition}
\newtheorem{rem}[thm]{Remark}

\def\sg{\sigma}

\begin{document}
\pagestyle{myheadings}

\title{Statistics on Wreath Products\\ and Generalized Binomial-Stirling 
Numbers}

\author{Amitai Regev\thanks{Partially supported by Minerva
Grant No. 8441 and by EC's IHRP Programme, within the Research Training Network ``Algebraic Combinatorics
in Europe'', grant HPRN-CT-2001-00272}\\
Department of Mathematics\\
The Weizmann Institute of Science\\
Rehovot 76100, Israel\\
{\em regev@wisdom.weizmann.ac.il}\\
\and
Yuval Roichman\thanks{Partially supported by EC's IHRP Programme, within the Research Training Network ``Algebraic Combinatorics
in Europe'', grant HPRN-CT-2001-00272}\\
Department of Mathematics\\
Bar Ilan University\\
Ramat Gan 52900, Israel\\
{\em yuvalr@math.biu.ac.il}}
\date{Apr 18, 2004}

\maketitle

\bibliographystyle{acm}

\begin{abstract}
Various statistics on wreath products are defined via canonical
words, ``colored'' right to left minima and ``colored'' descents.
It is shown that refined counts with respect to 
these statistics have nice recurrence formulas of binomial-Stirling type.
These extended Stirling numbers determine 
(via matrix inversion) dual systems, 
which are also shown to have combinatorial 
realizations within the wreath product.
The above setting also gives rise
to MacMahon type equi-distribution theorem 
over subsets with prescribed statistics.
\end{abstract}

\section{Introduction}


\bigskip

This paper was motivated by~\cite{RR1,RR2}, 
where we studied a variety of natural 
statistics on the symmetric group $S_n$ which
 generalized the {\it length},  
the {\it major} and other statistics. In particular, new statistics 
based on canonical 
presentations by the Coxeter generators were introduced.
Then the various 
Stirling numbers were obtained as cardinalities of certain subsets
of $S_n$ defined via these statistics.
For example, the Stirling numbers of the second kind are cardinalities
of subsets of permutations with prescribed number of left-to-right
minima and descents.
Refinements of the classical MacMahon-type equi-distribution 
theorems~\cite{MM} -- in the spirit of the results of Foata-Sch\"utzenberger, 
Garsia-Gessel etc. -- were deduced. 

\medskip

In this paper the group of permutations $S_n$ is replaced  by the
wreath product $C_a\wr S_n$, whose elements are called 
{\it ``colored permutations''}. Here $C_a$ is the cyclic group with 
$a$ elements. We study canonical presentations in wreath products
and introduce statistics counting the number of ``long" and of ``short" factors
in these presentations. These numbers essentially count number of certain 
right to left minima in colored permutations. It is shown that enumeration
of elements in wreath products with respect to these 
(and to these and descent) statistics
have nice recurrence formulas of binomial-Stirling type.
In particular, we present a wreath product extension of Stirling numbers 
of first and second kinds~\cite{Stir}, 
interpret these numbers in the wreath product,
 and prove MacMahon type equi-distribution theorem 
over subsets with prescribed statistics.

\medskip

Fix four integers $a,d,r,\ell\in Z$ and let
$g(n,k)=g_{a,d,r,\ell}(n,k)$ be the numbers determined
by the following recurrence:

\medskip

$g(0,0)=1$ ~and
\begin{eqnarray}\label{eex1}
g(n,k)=(an+dk-r)\cdot g(n-1,k)+\ell\cdot g(n-1,k-1),
\end{eqnarray}
and $g(n,k)=0$ if $k<0$ or $n<k$.

\medskip

The numbers 
$g_{a,d,r,\ell}(n,k)$
combine and generalize the binomial
coefficients and the Stirling numbers, see Section~\ref{s0}.
For example, $g_{1,0,1,1}(n,k)$ are the {\it signless Stirling numbers 
of the first kind}, $g_{0,1,0,1}(n,k)$ are the {\it  Stirling numbers 
of the second kind}, and  $g_{0,0,-1,1}(n,k)$ are the {\it  binomial 
coefficients}.
One of the main goals of this paper is to realize the numbers 
$g_{a,d,r,\ell}(n,k)$ 
via statistics on the wreath-products
$C_a\wr S_n$.


\bigskip

%


For a positive integer $a$ and a subset $L\subseteq \{0,\dots,a-1\}$
of cardinality $\ell$ let
\[
A_{L}(n,k):=\{\sg\in C_a\wr S_n\mid \m _L (\sg)=k\},
\]
and
\[
B_L(n,k):=\{\sg\in  C_a\wr S_n\mid des_L(\sg)=\m_L(\sg)=k\},
\]
where $\m_L(\sg)$ is the number of $L$-colored right to left minima, see Definition~\ref{uu1}.2, and $des_L(\sg)$
is the number of descents with respect to the $L$-order,
see Definitions~\ref{ord} and \ref{s11}. Then

\begin{thm}\label{main0} 
(See Corollary~\ref{ch2} and Theorem~\ref{sec8})
$$
g_{a,0,\ell,\ell}=\# A_L(n,k)
$$
and
$$
g_{0,a,\ell-a,\ell}(n,k)=\# B_L(n,k)
$$ 
\end{thm}

These two systems are essentially dual. This is

\begin{thm}\label{main1} (See Theorem~\ref{dual})
For every positive integers $a$, $N$, and every subset 
$L\subseteq \{0,\dots,a-1\}$ of size $\ell$, 
let $s_{L,N}$ 
be the $N\times N$ matrix
whose entries are given by
$$
s_{L,N}(n,k):=
{(-1)^{n-k}
\over \ell^n}
\cdot \# A_L(n,k)
\qquad (0\le k,n\le N)
$$

and $S_{L,N}$ be the $N\times N$ matrix whose entries are
$$
S_{L,N}(n,k):=
{1\over \ell^n}\cdot 
\# B_L(n,k)
\qquad (0\le k,n\le N).
$$
Then 
$$
S_{L,N}^{-1}= s_{L,N}.
$$
\end{thm}


To prove this theorem we apply a general decomposition and inversion 
theorems for linear recurrences, see Theorems~\ref{pt1} and~\ref{c1} below.
These theorems are closely related to results of
Milne and followers~\cite{Mil, MB, N}.

\bigskip


In Section~\ref{eqd} we apply the above setting 
to show that the length function and the
flag major index are equi-distribution over subsets of $B_n=C_2\wr S_n$
with prescribed colored right-to-left minima, see 
Corollaries~\ref{equid0},~\ref{equid2}.
This result is a type $B$-analogue of a recent theorem of Foata and 
Han for the symmetric group \cite[(1.5)]{FH} and refines a recent 
result of Haglund, Loehr and Remmel \cite[Theorem 4.5]{HLR}.

\bigskip

The rest of the paper is organized as follows.

Basic facts about wreath products are given in Sections~\ref{prel} and~\ref{wrp}.
In Section~\ref{stats} statistics on $C_a\wr S_n$ based on  
canonical words and on ``colored'' orders are introduced . 
Generalized Stirling numbers are interpreted combinatorially in 
Sections~\ref{col5},~\ref{col6} and~\ref{col9} , and are formally studied in 
Section~\ref{s0} and Appendix 2. The main equi-distribution theorem, 
Theorem~\ref{equid01}, is given in Section~\ref{eqd}.

\section{Preliminaries}\label{prel}

\noindent{\bf The wreath product $C_a\wr S_n$.}

Let $G$ be a group. Recall that the elements of the wreath product
$G\wr S_n$ are of the form $\sg= (( x_1,\ldots,x_n),p)$ where $x_i\in G$ and
$p\in S_n$; multiplication is given by 
\[
(( x_1,\ldots,x_n),p)(( y_1,\ldots,y_n),q)=
((x_1y_{p^{-1}(1)},\ldots,x_ny_{p^{-1}(n)}),pq).
\]
Let $A$ be the set 
$A:= G\times \{1,\ldots,n\}\equiv\{xj\mid x\in G,\ 1\le j\le n\}$.
We identify $(( x_1,\ldots,x_n),p)$ with the function 
$(( x_1,\ldots,x_n),p)\equiv f:A\to A$, given by 
\[
f:t j\to (t x_{p (j)})p (j)
\]
for all $t\in G$ and $1\le j\le n$.
When $G$ is Abelian one verifies easily that if, also, 
$g\equiv (( y_1,\ldots,y_n),q)$ then 
$fg\equiv (( x_1,\ldots,x_n),p)(( y_1,\ldots,y_n),q)$. This 
justifies the above identification. 
We therefore represent the element $\sg =(( x_1,\ldots,x_n),p)\in G\wr S_n$
by the $n$-tuple $[x_{p (1)}p(1), \ldots, x_{p (n)}p(n) ]$:
\[
\sg =(( x_1,\ldots,x_n),p)\equiv
[x_{p (1)}p(1), \ldots, x_{p (n)}p(n) ]=[\sg(1),\ldots,\sg(n)],
\] 
and we denote $|\sg |:=p$. Note that if 
$\sg =[y_1j_1,\ldots,y_nj_n]$ where $y_i\in G$ and $1\le j_i\le n$, 
then $| \sg|=p=[j_1,\ldots,j_n]$. 
Let $z_i\in G$, $p\in S_n$ and let  $\sg=[z_1p(1),\ldots,z_np(n)]$, then 
$\sg^{-1}=[z_{p^{-1}(1)}^{-1}p^{-1}(1),\ldots,z_{p^{-1}(n)}^{-1}p^{-1}(n)]$.

\medskip

In this paper we consider the wreath products $C_a\wr S_n$, where $C_a$
is the (multiplicative) cyclic group of order $a$: 
$\alpha:=e^{2\pi i \over a}$, and 
\[
C_a:=\{\alpha^t \mid 0\le t\le a-1\}.
\]
The elements of $C_a\wr S_n$ are identified with 
``$a$--colored'' permutations,
namely those permutations
$\sg$ of the set $A=C_a\times\{1,\ldots,n\}$ satisfying
\[
\sg (\beta j)=\beta \sg (j),\quad \beta\in C_a,\quad 1\le j\le n.
\]
We write $\sg=[\sg(1),\ldots,\sg(n)]$. For each $j$, 
$\sg(j)=\alpha^{t_j}\cdot |\sg(j)|$ and has {\it color} $\alpha^{t_j}$; it is ``colorless'' if $t_j=0$.

\medskip

\noindent
{\bf Cycle decomposition}:
Let $\sg =((x_1,\ldots,x_n),p)\in C_a\wr S_n$. 
The cycle decomposition of $p$ induces the corresponding 
decomposition of $\sg$: 
If $p=p_1\cdots p_r$ is 
the cycle decomposition of $p$, 
and $p_i=(b^{(i)}_{1},\ldots,b^{(i)}_{m})$ (in the ordinary cycle notation for $S_n$),
for each $1\le i\le r$ let
\[
y_j^{(i)}:=\cases
{x_j, &  if $j\in\{b_{1}^{(i)},\ldots,b_{m}^{(i)}\}$; \cr
1, & otherwise. \cr} 
\]
Then $\sg^{(i)}:=((y_1^{(i)},\ldots,y_n^{(i)}), p_i)$ are the corresponding cycles of $\sg$, and $\sg=\sg^{(1)}\cdots\sg^{(r)}$ 
is the {\it cycle decomposition} of $\sg$ 
The product
$x_{b^{(i)}_1}\cdots x_{b^{(i)}_m}\in C_a$ is uniquely determined (since $C_a$
is Abelian), and is called the {\it color of that cycle} 
of $p=|\sg |$.

%

\bigskip

\noindent
{\bf Generators and length.}
Let $s_i=(i,i+1)\in S_n$, $i=1,\ldots n-1$, denote the Coxeter generators
of $S_n\subset C_a\wr S_n$. In addition, $s_0\in C_a\wr S_n$
is the element given by
\[
s_0=((\alpha,1,\ldots,1),1)\equiv  
[\alpha, 2,3,\ldots ,n]
\]
where  $\alpha=e^{2\pi i\over a}$.\\


The following easy fact is well known.
\begin{fac}\label{b1}
Let $\sg =[b_1,\ldots,b_n]\in C_a\wr S_n$. 
\begin{enumerate}
\item
$\sg s_0=[\alpha b_1,b_2,\ldots,b_n]$.
\item
Let $1\le i\le n-1$, then
$\sg s_i=[b_1,\ldots,b_{i+1},b_i,\ldots,b_n]$
\end{enumerate}

\end{fac}
The set
$S=\{s_0,s_1,\ldots,s_{n-1}\}\subseteq C_a\wr S_n$ 
generates $C_a\wr S_n $ 
(this follows, for example, from
Proposition~\ref{canonical}).

The {\it length} of an element $\sg\in C_a\wr S_n$, denoted $\ell(\sg)$,
is the minimum length of an expression of $\sg$ as a product
of elements in the above generating set $S$.

\section{Canonical Presentation in Wreath Products}\label{wrp}

Consider the following subsets of elements in $C_a\wr S_n$. First, let
$R_0=C_a$.
Given $1\le j\le n-1$, let
$$
R_j^0:=\{1, s_j, s_j s_{j-1}, \dots, s_j\cdots s_1\}.
$$
For $1\le t\le a-1$ let
$$
R_j^t:=\{s_j\cdots s_0^t,\; s_j\cdots s_0^t s_1,\; s_j\cdots s_0^t s_1 s_2,\;
 \dots,\;
s_j\cdots s_0^t s_1 \cdots s_j\}
$$
and
$$
R_j:=\bigcup_{t=0}^{a-1} R_j^t.
$$
Note  that $|R_j|=a\cdot (j+1)$, hence

\[
\prod_{j=0}^{n-1}|R_j|=a^n\cdot n!=|C_a\wr S_n|.
\]
\begin{pro}\label{canonical}
Every element $\sigma\in C_a\wr S_n$ has a unique presentation
$$
\sigma=w_0\cdots w_{n-1}
$$
where, for every $0\le j\le n-1$, $w_j\in R_j$.
\end{pro}

\noindent{\bf Proof.}
By induction on $n$ and by Fact~\ref{b1}.
Recall that every element $\sg\in C_a\wr S_n$ may be interpreted as a 
colored permutation
$[\sg(1),\dots,\sg(n)]$. It follows from this interpretation that
every element $\sg\in C_a\wr S_n$ is obtained in a unique way by 
inserting colored $n$ (namely $e^{2\pi i t\over a}n$
for some $0\le t\le a-1$)
into a colored permutation $\bar\sg \in C_a\wr S_{n-1}$.
Now, if $\sg(j)=e^{2\pi it\over a} n$ and $t=0$ then $\sg=\bar\sg w_{n-1}$,
where
$$
w_{n-1}=\cases
{1, &  if $j=n$; \cr
s_{n-1}\cdots s_{j}, & if $j<n$. \cr} \in R_{n-1}^0.
$$

\noindent
If $\sg(j)=e^{2\pi it\over a} n$ and $0< t\le a-1$
then $\sg=\bar\sg w_{n-1}$
where
$$
w_{n-1}=s_{n-1}\cdots s_0^t \cdots s_{j-1}\in R_{n-1}^t.
$$
This proves ``existence''. Uniqueness now follows by a standard counting
argument.
\qed

\medskip

\begin{df}\label{ca1}
Call the above presentation $\sg=w_0\cdots w_{n-1}$
in Proposition~\ref{canonical} the canonical presentation --
or the canonical word -- of $\sg=w_0\cdots w_{n-1}$.
\end{df}
%


\begin{pro}\label{Bagno}
Write the above canonical word explicitly:
$\sg=w_0\cdots w_{n-1}=s_{i_1}\cdots s_{i_r}$, then $r$ 
is the minimum length of an expression of $\sg$ as a product
of elements in  $S=\{s_0,s_1,\ldots,s_{n-1}\}$,
i.e.~the length of $\sg$ is $\;\ell (\sg ) = r$. 
\end{pro}

For a proof see e.g.~\cite[Ch. 3.3]{B}.

\begin{cor}\label{ca2}
Let $\sg=w_0\cdots w_{n-1}$ be the canonical word of $\sg\in C_a\wr S_n$,
then $\ell (\sg)=\ell(w_0)+\cdots +\ell(w_{n-1})$. In particular, if
$\bar\sg\in C_a\wr S_{n-1}$ and $r\in R_{n-1}$ then
$\ell(\bar\sg r)=\ell(\bar\sg)+\ell(r)$.
\end{cor}

\section{Statistics on Colored Permutations}\label{col0}
\label{stats}

In this section we 
introduce
various statistics on $C_a\wr S_n$ based
on canonical words, on right-to-left-minima, and on
certain descent sets $Des_L$.

\subsection{Right to Left Minima}\label{Right}


Recall from Section~\ref{prel}
the notation $\al:=e^{2\pi i\over a}$
and  $|\sigma|$ (for every $\sg\in C_a\wr S_n$).


%

%
%
\begin{df}\label{uu1}
\begin{enumerate}
\item
Let $p=[j_1,\ldots,j_n]\in S_n$. 
Define ${\stackrel{\longleftarrow}{\mbox{Min}}}(p)\subseteq\{1,\ldots,n\}$ as follows: 
\[
{\stackrel{\longleftarrow}{\mbox{Min}}}(p)=\{j_i\mid j_i\;\mbox{is a 
r.t.l.min 
in $[j_1,\ldots,j_n]$}\}.
\]
Here and on r.t.l.min stands for right to left minimum.
\item
Let  $L\subseteq \{1,\dots, a-1\}$. Let $\sg\in C_a\wr S_n$ be a colored
permutation, and write $\sg=[b_1,\ldots,b_n]$.  
Define ${\stackrel{\longleftarrow}{\mbox{Min}}}_L(\sg)\subseteq\{1,\ldots,n\}$ 
as follows:
\[
{\stackrel{\longleftarrow}{\mbox{Min}}}_L(\sg)=\{|b_i|\mid |b_i|\;\mbox{is a r.t.l.min in $|\sg|$,
and $b_i=\al^u|b_i|$ 
for some $u\in L$}\}.
\]
Finally  denote ${\stackrel{\longleftarrow}{\mbox{min}}}_L(\sg)=|{\stackrel{\longleftarrow}{\mbox{Min}}}_L(\sg)|$.
\end{enumerate}
\end{df}
For example let $\sg=[\al 3, \al ^3 5, 1,\al^2 2 ,\al 4]\in C_4\wr S_5$, then $|\sg|=[3,5,1,2,4]$ and 
${\stackrel{\longleftarrow}{\mbox{Min}}}_{\{0,1,2,3\}}(\sg)=\{1,2,4\}$,
${\stackrel{\longleftarrow}{\mbox{Min}}}_{\{1,2\}}(\sg)=\{2,4\}$,
 ${\stackrel{\longleftarrow}{\mbox{Min}}}_{\{0,3\}}(\sg)=\{1\}$, and ${\stackrel{\longleftarrow}{\mbox{Min}}}_{\{0,1,2\}}(\sg)=\{1,2,4\}$.

\begin{pro}\label{uu2} 
\begin{enumerate}
\item
Let $p=[j_1,\ldots,j_n]\in S_n$, let $v_0=1$ and let
$p=v_0v_1\cdots v_{n-1}$ be its canonical presentation. Then 
$j_i$ is a r.t.l.min in $[j_1,\ldots,j_n]$ if and only if $v_{j_i-1}=1$.
\item
Let  $\sg\in C_a\wr S_n$ and
$\sg=w_0\cdots w_{n-1}$   $\;(\forall i\  w_i\in R_i)$ be its canonical word.
Also let  
$|\sg|=v_1\cdots v_{n-1}$ be the canonical presentation of $|\sg|$. 
For each $0\le u\le a-1$ and $1 \le j\le n-1$ denote
$$
r_{u,j}:=\cases
{s_j\cdots s_0^u\cdots s_j\in R_j &  if $u\ne 0$ \cr
1 & if $u=0$. \cr}
$$
Then $v_i=1$ if and only if $w_i=r_{u,i}$ for some $u$.
\item
Let $L\subseteq \{0,\dots, a-1\}$. Then
$$
{\stackrel{\longleftarrow}{\mbox{Min}}}_L(\sg)=\{0\le i \le n-1| \exists u\in L,\ w_i  = r_{u,i}\}.
$$
\end{enumerate}
\end{pro}

\noindent{\bf Proof -} is standard
(by induction on $n$) and is left to the reader.

\begin{cor}\label{uu4}
Let $\bar\sg\in C_a\wr S_{n-1}$ and $r=w_{n-1}\in C_a\wr S_{n-1}$ (hence
$\sg\in C_a\wr S_n$). 
Let 
\begin{equation}\label{df-kl}
K_L(r)=K_L(w_{n-1}):=\cases
{\{n-1\} &  if ~$\exists ~u\in L~~w_{n-1}=r_{u,n-1}$ \cr
\phi & otherwise. \cr}
\end{equation}
Then 
\begin{equation}\label{uu41}
\M_L(\sg)=\M_L(\bar\sg )\cup K_L(w_{n-1}),~~\mbox{a disjoint union},
\end{equation}
\end{cor}

\bigskip


%

%

\subsection{The Order $<_L$ and the $L$-Descent Set}

Notice that $C_a\wr S_n$ is identified with the permutations
$\sg$ of the set\\
$\{\al^vj\mid 0\le v\le a-1,\;1\le j\le n\}\cup\{0\}$, where by definition,
$\sg (0)=0$, and $\sg(\al^vj)=\al^v\sg(j)$.

\begin{df}\label{ord}
A subset $L\subseteq\{0,\ldots,a-1\}$ determines a linear order $<_L$
on  $\{\al^vj\mid 0\le v\le a-1,\; 0\le j\le n\}\cup\{0\}$ 
as follows:

Let $U=\{0,\ldots,a-1\}\setminus L$ be the complement of $L$ in
$\{0,\ldots,a-1\}$.

If $v\in L$ then $\al^v j<_L 0$ for every $1\le j\le n$.
If $v\in U$ then $\al^v j>_L 0$ for every $1\le j\le n$.

For $v,u\in L$ (not necessarily distinct)
and $i\ne j\in\{1,\dots,n\}$,\\
$\al^v i<_L \al^u j$ if and only if $i > j$ (``reverse order'').

For $v,u\in U$ (not necessarily distinct)
and $i\ne j\in\{1,\dots,n\}$,\\
$\al^v i<_L \al^u j$ if and only if $i < j$.

Then, for each $1\le j\le n$,
order each subset $\{\al^vj\mid v\in L\}$ 
(and  each subset $\{\al^vj\mid v\in U\}$)
in an arbitrary linear order.

%
This yields a linear order $<_L$ 
on the set
$\{\al^vj\mid 0\le v\le a-1,\; 0\le j\le n\}\cup\{0\}$.
\end{df}

\medskip

For example let $a=4$ and $L=\{2,3\}$, then $U=\{0,1\}$. We can choose 
the following order
\[
\al^2n<_L\al^3n<_L\al^2(n-1)<_L\al^3(n-1)<_L\cdots<_L\al^2<_L\al^3<_L0~~~~~~~~~~~~~~~~
\]
\[
~~~~~~~~~~~~~~~~~~~
0<_L\al<_L1<_L\al2<_L2<_L\cdots <_L\al(n-1)<_L(n-1)<_L\al n<_Ln.
\]
The following is an obvious property of this order.

\begin{fac}\label{u10}
Let $\bar\sg =[\bar\sg(1),\ldots,\bar\sg(n-1)]
\in C_a\wr S_{n-1}$ 
and let $0\le v\le a-1$.\\
If $v\in L$ then $\al^vn <_L\bar\sg(1),\ldots, \bar\sg(n-1)$;\\
if $v\in U$ then $\al^vn >_L\bar\sg(1),\ldots, \bar\sg(n-1)$.
\end{fac}


\begin{df}\label{s11}
The $L$-descent set of $\sg\in C_a\wr S_n$ is
\[
Des_L(\sg):=\{0\le i\le n-1\;\mid \;\sg (i)>_L\sg(i+1)\}.
\]
The $L$-descent number is 
$$
des_L(\sg):=|Des_L(\sg)|.
$$
If $L$ consists of one element $u\in \{0,\dots, a-1\}$ then
we denote $<_u$, $\D_u$, $des_u$.

\end{df}
%





The following notion is the natural
analogue of the standard descent sets of Weyl and Coxeter groups.

\begin{df}\label{sd1}
For $\sg\in C_a\wr S_n$ let the {\it standard descent set} be
$$
\D (\sg):=\{ 0\le i\le n-1|\ \ell(\sg s_i)<\ell(\sg)\}.
$$
\end{df}


It should be noted that
the $u$-descent set, $\D_u$, defined above,
may also be interpreted via the generators.

\begin{pro}\label{ddes-L}
For every $\sg\in C_a\wr S_n$ and every $0\le u\le a-1$ 
$$
\D_{u}(\sg) = 
\{
0\le i\le n-1|\ 
\ell(v_u^{-1}\sg s_i)>
\ell(v_u^{-1}\sg)
\}, 
$$
where $v_u:= ((\al^u,\dots,\al^u),id)=[\al^u 1, \al^u 2,\dots,  \al^u n]$.
\end{pro}

\noindent
{\bf Proof -} 
is given in an appendix (Section~\ref{app1}).

\smallskip

\begin{exa}\label{ddes-L-example}
\end{exa}

\noindent
{\bf 1)} $L=\{0,\dots,a-1\}$. By definition,
$$
\D_{\{0,\dots,a-1\}}(\sg)=\D (|\sg|)=\{0\le i\le n-1|\ |\sg(i)|>|\sg(i+1)|\} 
$$
the standard descent set of $|\sg|$.

\noindent
{\bf 2)} $L=\emptyset$.  $\D_\emptyset(\sg)$ is the complement of the standard
descent set of $|\sg|$ (the ascent set of $|\sg|$).

\noindent
{\bf 3)} $L=\{1,\dots,a-1\}$. 
By Proposition~\ref{ddes-L}, since $v_0$ is the identity element 
$$
\D_{\{1,\dots,a-1\}}(\sg)=
\{0,\dots,n-1\}\setminus \D_0(\sg)=\{0\le i\le n-1|\ \ell(\sg s_i)<\ell(\sg)\}
$$
the standard descent set of $\sg$.

\noindent
{\bf 4)} $L=\{0,\dots,a-2\}$. $v_{a-1}$ is the longest element 
in $C_a\wr S_n$ and $\D_{\{0,\dots,a-2\}}$ is the complement
of the standard descent set; namely, the ascent set of $\sg$.

\bigskip

\begin{lem}\label{p2}
Let $L\subseteq\{0,\ldots,a-1\}$ then,
for any $\sg\in C_a\wr S_n$, 
$$
des_L(\sg)\ge\ \der_L(\sg).
$$
\end{lem}
{\bf Proof.} Let $\Der_L(\sg)=\{i_1,\ldots,i_k\}$, and show that
for each $1\le j\le k-1$, $\sg(i_j)>_L\sg(i_{j+1})$. Indeed, each
$\sg (i_t)=\al^{v_t}|\sg (i_t)|$, $v_t\in L$, and $|\sg (i_t)|$ is
a r.t.l.min of $|\sg|$. Therefore $|\sg (i_j)|<|\sg (i_{j+1})|$,
so
\[
\sg(i_j)=\al^{v_j}|\sg (i_j)|>_L\al^{v_{j+1}}|\sg (i_{j+1})|=\sg
(i_{j+1}),
\]
as was claimed. By the transitivity of the linear order $>_L$,
there must be an $L$-descent of $\sg$ between these two indices
$i_j$ and $i_{j+1}$. This contributes (at least) $k-1$
$L$-descents to $Des_L(\sg)$. By definition,
$\sg(0)=0>_L\sg(i_1)$, and this contributes at least one more
$L$-descent of $\sg$. \qed

\medskip

{\bf Note} that here we have to allow $0\in Des_L(\sg)$.

%
%
%
%
%

\section{``Colored" Stirling Numbers of the First Kind}\label{col5}

In this section we point on connections between statistics on
colored permutations, defined above, and certain generalized Stirling numbers of the first kind.

\begin{pro}\label{ch1}
Let $L\subseteq \{0,\ldots,a-1\}$, $|L|=r$. Then
\[
\sum_{\sg\in C_a\wr S_n} q^{\m _L (\sg)} =(rq+a-r)(rq+2a-r)\cdots (rq+na-r).
\]
\end{pro}
{\bf Proof.} By Corollary~\ref{uu4} it suffices to show that for every $n$ 
\[
\sum_{w_{n-1}\in R_{n-1}} q^{|K_L (w_{n-1})|}=rq+na-r.
\]
Indeed, by definition
(\ref{df-kl}) (in Corollary~\ref{uu4})
\[
\sum _{w_{n-1}\in R_{n-1}} q^{|K_L (w_{n-1})|}=rq+|R_{n-1}|-r=rq+na-r.
\]
\qed

\begin{cor}\label{ch2}
Let $r=|L|$ as above, and denote
\[
g_L(n,k):=\#\{\sg\in C_a\wr S_n\mid \m _L (\sg)=k\}.
\]
Then $g_L(n,k)$ satisfies the following recurrence:
\[
g_L(n,k)=(an-r)\cdot g_L(n-1,k)+r\cdot g_L(n-1,k-1).
\]
Thus, by Equation~(\ref{eex1}), $g_L(n,k)=g_{a,0,r,r}(n,k)$, so
\[
g_{a,0,r,r}(n,k)=\#\{\sg\in C_a\wr S_n\mid \m _L (\sg)=k\}.
\]
\end{cor}

{\bf Proof.} By Proposition~\ref{ch1}
\[
\sum_k g_L(n,k)q^k=
\sum_{\sg\in C_a\wr S_n} q^{\m _L (\sg)} =(rq+a-r)(rq+2a-r)\cdots (rq+na-r).
\]
Thus
\[
\sum_k g_L(n,k)q^k=(rq+na-r)\sum_k g_L(n-1,k)q^k=~~~~~~~~~~~~~~~~~~~
\]
\[
~~~~~~~~~~~~~~~~~~~
=(na-r)\sum_k g_L(n-1,k)q^k+\sum_k r\cdot g_L(n-1,k-1)q^k,
\]
and the proof follows. \qed

\medskip

{\bf Note} that when $a=|L|=r=1$, $g_L(n,k)$ are the signless Stirling numbers
of the first kind. In Section~\ref{s0} we study similar but more general 
such recurrences. 

\medskip

Recall from Section 2 that the cycles of $\sg\in C_a\wr S_n$ are ``colored''
by elements of $C_a$.
\begin{df}\label{ch3}
Given $L\subseteq\{1,\ldots,n\}$ and $\sg\in C_a\wr S_n$, we say that  a
cycle of $\sg$ is $L$--colored if its color belongs to $L$.
\end{df}
%


\begin{cor}\label{ch4}
The number of
elements $\sg\in C_a\wr S_n$ with exactly $k$ ~r.t.l.min 
of $| \sg|$ which are $L$--colored, $g_L(n,k)$, is also the number of
elements $\sg\in C_a\wr S_n$ with exactly $k$ cycles which are 
$L$--colored.
\end{cor}
\noindent{\bf Proof.} 
The proof is a natural extension of~\cite[p. 17]{ECI}.
The following notion will be used in the proof.
Let $\sg=((x_1,\ldots,x_n),p)\in C_a\wr S_n$, and let
$\gamma=(b_1,\ldots, b_m)$ be a cycle of $p=|\sg|$. 
Assume w.l.o.g.~that the last element $b_m$ is minimal, then
the color of $b_m$, 
$x_{b_m}$, will be called 
{\it the right--color} of the cycle $\gamma$. 
A cycle is {\it right $L$--colored}, for $L\subseteq \{0,\dots,a-1\}$,
if its right--color belongs to $\{\alpha^u|\ u\in L\}$.

Let $G'_L(n,k)$ denote the set of 
elements $\sg\in C_a\wr S_n$ with exactly $k$ cycles which are right
$L$--colored and
\[
G_L(n,k):=\{\sg\in C_a\wr S_n\mid \m _L (\sg)=k\}.
\]
We first construct a bijection 
\[
G'_L(n,k) \longleftrightarrow G_L(n,k).
\]
Given $\sg'=((x_1,\ldots,x_n),p')\in G'_L(n,k)$, reorder the cycles 
in $p'=|\sg'|$ such that each cycle in $|\sg'|$ is
written with its smallest element last (i.e.~rightmost), and the 
cycles are written in increasing order of their smallest element.
By assumption, exactly $k$ of these smallest elements are $L$--colored. 
Let $p$ be the permutation obtained from $p'$ by erasing the parenthesis
of the cycles, and let $\sg=((x_1,\ldots,x_n),p)$. Clearly, in $p=|\sg|$, 
those smallest elements are now r.t.l.min, and in $\sg$ they have 
the same colors as in 
$\sg'$, namely  exactly $k$ of these r.t.l.min are $L$--colored. Thus
$\sg\in G_L(n,k)$. That correspondence can be reversed by parenthesizing 
$p\in S_n$ according to its r.t.l.min, 
therefore the above is a bijection.

Let $G''_L(n,k)$ denote the set of 
elements $\sg\in C_a\wr S_n$ with exactly $k$ cycles which are
$L$--colored. There is a rather obvious bijection 
\[
G''_L(n,k) \longleftrightarrow G'_L(n,k)
\]
as follows. Given $\sg''=((x_1,\ldots,x_n),p'')\in G''_L(n,k)$,
let $(b_1,\ldots,b_m)$ be a cycle of $p''$ with $b_m$ minimal,
then replace $x_{b_m}$ by $x_{b_1}\cdots x_{b_m}$. Do it to each cycle.
This clearly maps $G''_L(n,k) \longrightarrow G'_L(n,k)$, with an 
obvious inverse map. This completes the proof.
\qed

\section{``Colored" Stirling Numbers of the Second Kind}\label{col6}

In this section 
we prove the second part of Theorem~\ref{main0} (Theorem~\ref{sec8} below).

\medskip

Throughout this section we assume that 
$L\subseteq\{0,1,\ldots,a-1\}$, with the
corresponding linear order $<_L$ 
as above.
%
%
%
%

\begin{lem}\label{sec3}
Let $\sg=w_0\cdots w_{n-1}$ (canonical presentation),
$\bar\sg=w_0\cdots w_{n-2}$, so $\sg= \bar\sg w_{n-1}$. Then
$des_L(\sg)\ge_L des_L (\bar\sg)$.
\end{lem}

\noindent{\bf Proof.} Recall that $\sg$ is obtained from $\bar\sg$ by
inserting some $\al^vn$ into $\bar\sg$. Thus, for certain
$b_1,\ldots,b_{n-1}\in C_a \cdot\{1,\ldots,n-1\}$ and $1\le t\le
n-1$,
\[
\bar\sg=[b_1,\ldots,b_{n-1},n] \quad\mbox{and}\quad
\sg=[b_1,\ldots,b_t,\al^vn,b_{t+1},\ldots,b_{n-1}].
\]
Since the $L$-order is linear, if $b_t>_Lb_{t+1}$ then either
$b_t>_L\al^vn$ or (and/or) $\al^vn>_Lb_{t+1}$, which implies the
proof. \qed
\begin{lem}\label{sec4}
With the notation of the previous Lemma,
\begin{enumerate}
\item
if $\sg(n)=\al^vn$ for some $v\in L$ then
$\M_L(\sg)=\M_L(\bar\sg)\cup\{n\}$, hence\\
$\m_L(\sg)=\m_L(\bar\sg)+1$;
\item
if $\sg(n)\ne \al^vn$ for any $v\in L$ then
$\m_L(\sg)=\m_L(\bar\sg)$.
\end{enumerate}
\end{lem}
\noindent{\bf Proof.} The lemma is an immediate consequence of Corollary~\ref{uu4}.

\qed

The following is a key observation here.

%

\begin{lem}\label{sec5}
Let $\sg=\bar\sg w_{n-1}$ as above, and assume\\
$des_L(\sg)=\m_L(\sg)=k$.
\begin{enumerate}
\item
If $\sg (n)=\al^vn$ for some $v\in L$ then
$des_L(\bar\sg)=\m_L(\bar\sg)=k-1$.
\item
If $\sg (n)\ne\al^vn$ for any $v\in L$ then
$des_L(\bar\sg)=\m_L(\bar\sg)=k$.
\end{enumerate}
\end{lem}
\noindent{\bf Proof.}

\noindent
1. Assume $\sg(n)=\al^vn$, $v\in L$. By 
Lemma~\ref{sec4}.1,
$\m_L(\bar\sg)=\m_L(\sg)-1$. Clearly in that case
$Des_L(\sg)=Des_L(\bar\sg)\cup\{n-1\}$, hence also
$\des_L(\bar\sg)=k-1$.

\medskip

\noindent
2. If $\sg (n)\ne \al^vn$ for any $v\in L$ then, by Lemma~\ref{sec4}.2,
$\m_L(\sg)= \m_L(\bar\sg)$.
Therefore by Lemmas~\ref{sec3} and~\ref{p2},
\[
k=\des_L(\sg)\ge des_L(\bar\sg)\ge \m_L(\bar\sg)=\m_L(\sg)=k,
\]
forcing equality. Thus $des_L(\bar\sg)=\m_L(\bar\sg)=k$. \qed

\medskip

%
%

\begin{lem}\label{sec7}
Let $\bar\sg\in C_a\wr S_{n-1}$,
$\bar\sg=[\bar\sg(1),\ldots,\bar\sg(n-1)]$. 
\\
Assume that $des_L(\bar\sg)=\m_L(\bar\sg)=k$ and let
$Des_L(\bar\sg)=\{i_1,\ldots,i_k\}$.
\begin{enumerate}
\item
If $v\not\in L$ then there are exactly $k+1$ elements 
$\sg\in C_a\wr S_n$, such that $\sg=\bar\sg w_{n-1}$ for some 
$w_{n-1}\in R_{n-1}$ and $des_L(\sg)=\m_L(\sg)=k$.

\item
If $v\in L$ then 
 there are exactly $k$ such $\sg$'s ``over''
$\bar\sg$ satisfying\\ $des_L(\sg)=\m_L(\sg)=k$.
\end{enumerate}
\end{lem}

\noindent{\bf Proof.}
Fix some $b_n=\al^vn$ and insert it into  $\bar\sg$
to obtain $\sg=\bar\sg w_{n-1}$.

\noindent 1. $v\not\in L$, hence $\bar\sg(1),\ldots,\bar\sg(n-1)<_L b_n$. If
$b_n$ is inserted immediately to the right of some $\bar\sg(i_t)$
($1\le t\le k$) or in the last ($n$-th) position, then
$des_L(\sg)=k$. 
Also, by Lemma~\ref{sec4}.2, $\m_L(\sg)=\m_L(\bar\sg)=k$. Conversely, if
$des_L(\sg)=k$ then $b_n$ was inserted into one of these $k+1$
positions.

\noindent 2. $v\in L$,  hence $b_n<_L,0,
\bar\sg(1),\ldots,\bar\sg(n-1)$. If $b_n$ is inserted immediately
to the right of some $\bar\sg(i_t)$ ($1\le t\le k$)
then $des_L(\sg)=k$. 
Also, in this case $b_n$ is not inserted in the last position; 
by Lemma~\ref{sec4}.2, $\m_L(\sg)=\m_L(\bar\sg)=k$.
Conversely, if $des_L(\sg)=k$ then $b_n$ was
inserted into one of these $k$ positions.

Note that if $i_1\ne 0$ and $b_n$ is inserted in the first (left)
position then $0$ is an additional $u$-descent of $\sg$, since
$0>_Lb_n$.

\qed

%
%
%
%

\medskip

\begin{df}\label{ws2}
Let $f_L(0,0)=1$ and define
\[
f_L(n,k)=\#\{\sg\in  C_a\wr S_n\mid des_L(\sg)=\m_L(\sg)=k\}.
\]
\end{df}
\begin{thm}\label{sec8}
Let 
$\ell =|L|$. Then $f_L(n,k)$
satisfies the following recurrence:
\[
f_L(n,k)=(ak+a-\ell)\cdot  f_L(n-1,k)+\ell\cdot f_L(n-1,k-1).
\]
Thus $f_L(n,k)=g_{0,a,\ell-a,\ell}(n,k)$, so
\[
g_{0,a,\ell-a,\ell}(n,k)=\#\{\sg\in  C_a\wr S_n\mid des_L(\sg)=\m_L(\sg)=k\}.
\]
\end{thm}

\noindent{\bf Proof.} 
Let 
$$
B_L(n,k):=\{\sg\in  C_a\wr S_n\mid des_L(\sg)=\m_L(\sg)=k\},
$$ 
so $f_L(n,k)=\#B_L(n,k)$,
$$
C_L(n,k):=\{\sg=\bar\sg w_{n-1}\in B_L(n,k)\mid \bar\sg\in B_L(n-1,k-1)\},
$$
and
$$
D_L(n,k):=\{\sg=\bar\sg w_{n-1}\in B_L(n,k)\mid \bar\sg\in B_L(n-1,k)\}.
$$
By definition, $C_L(n,k)\cap D_L(n,k)=\emptyset$. 
By Lemma~\ref{sec5},
$$
B_L(n,k)=C_L(n,k)\cup D_L(n,k).
$$

The proof will follow, once we show that 
\begin{enumerate}
\item
$|C_L(n,k)|=\ell\cdot|B_L(n-1,k-1)|$ and
\item
$|D_L(n,k)|=(ak+a-\ell)\cdot  |B_L(n-1,k)|$,
\end{enumerate}

\smallskip

\noindent 1. By Lemma~\ref{sec5}, 
all elements in $C_L(n,k)$ which are obtained 
from an element $\bar\sg\in B_L(n-1,k-1)$
by inserting a colored $n$, are obtained by inserting
an $L$-colored $n$ at the last position: 
$\sg=[\bar\sg ,\al^vn]$, $v\in L$. This proves 1.

\smallskip

\noindent 2. Let $\bar\sg\in B_L(n-1,k)$:
$\;des_L(\bar\sg)=\m_L(\bar\sg)=k$ and insert  $\al^v n$
into $\bar\sg$ to obtain a permutation $\sg=\bar\sg w_{n-1}\in B_L(n,k)$. 
If $v\not\in L$ then, by 
Lemma~\ref{sec7}.1, there are exactly $k+1$ 
such permutations $\sg\in B_L(n,k)$. Since there are $a-\ell$ such $v$'s, we get $(a-\ell)(k+1)$ ~$\sg$'s. Similarly, Lemma~\ref{sec7}.2 
implies $k$ such
$\sg$'s when $v\in L$, namely a total of $\ell k$  ~$\sg$'s. Together, this
yields exactly $(a-\ell)(k+1)+\ell k=ak+a-\ell$ 
~$\sg$'s in $B_L(n,k)$ ``over'' each $\bar\sg\in B_L(n-1,k)$.
This proves 2.

 \qed

\begin{rem}\label{classical}
Letting $a=\ell=1$, $f_L(n,k)$ are the classical Stirling numbers of the second kind.
\end{rem}


\section{Equi-distribution in $B_n=C_2\wr S_n$}\label{eqd}

In this section we study the case of $B_n=C_2\wr S_n$, namely $a=2$.
We prove here an equi-distribution theorem between the length parameter
$\ell(\sg)$ and the flag-major index, see Definition~\ref{df-fm} below.

Here 
$L\subseteq \{0,1\}$ determines $<_L$.
In the case $L=\{1\}$ the natural order is preserved, and it is reversed
when $L=\{0\}$:
\[
-n<_1-(n-1)<_1\cdots <_1 -1<_1 0<_1 1<_1\cdots<_1 n ~~~~\mbox{and} 
\]
\[
n<_0n-1<_0\cdots <_0 1<_0 0<_0 -1<_0\cdots<_0 -n.
\]
These orders define the corresponding $\M_0$ and the $\M_1$ sets,
see Definition~\ref{uu1}.
In this section we show that the 
length function and the
flag major index are equi-distribution over subsets of $B_n=C_2\wr S_n$
with prescribed  
$\M_0$ and $\M_{1}$ sets, see Corollary~\ref{equid0} below.
This result is a type $B$-analogue of a recent theorem of Foata 
and Han for the symmetric group~\cite[(1.5)]{FH} and refines a 
recent result of Haglund, Loehr and Remmel~\cite[Theorem 4.5]{HLR}.

\bigskip

\begin{thm}\label{equid00}
For every positive integer $n$
$$
\sum\limits_{\sigma\in B_n}
 \prod\limits_{i\in \M_0(\sg)}x_i 
\cdot \prod\limits_{i\in \M_{1}(\sg)}t_i 
\cdot q^{\ell(\sigma)}
=
$$
$$
(x_1+qt_1)(x_2+q+q^2+q^3t_2)\cdots (x_n+q+q^2+\cdots+q^{2n-1}t_n).
$$
\end{thm}

\noindent{\bf Proof.} 
By induction on $n$.
Obviously, theorem holds for $n=1$. 

\noindent
By Proposition~\ref{canonical}, the l.h.s.~equals 
\[
\sum\limits_{\bar\sigma\in B_{n-1}} 
\sum_{r\in R_{n-1}}q^{\ell(\bar\sigma r)}
\cdot \prod\limits_{i\in \M_0(\bar\sigma r)}x_i 
\cdot \prod\limits_{i\in \M_{1}(\bar\sigma r)}t_i =Q
\]
By  Remark~\ref{ca2}
and Corollary~\ref{uu4}, $Q$ equals
\[
\left [ \sum\limits_{\bar\sigma\in B_{n-1}} 
\prod\limits_{i\in \M_0(\bar\sigma )}x_i 
\cdot \prod\limits_{i\in \M_{1}(\bar\sigma )}t_i 
\cdot q^{\ell(\bar\sigma )}
\right ]
\cdot \left [\sum_{r\in R_{n-1}}
\prod\limits_{i\in K_0(r)}x_i
\cdot \prod\limits_{i\in K_1(r)}t_i 
\cdot q^{\ell(r)}
\right ].
\]

Thus, by induction, it suffices to show that 
\[
\left [\sum_{r\in R_{n-1}}
\prod\limits_{i\in K_0(r)}x_i
\cdot \prod\limits_{i\in K_1(r)}t_i 
\cdot q^{\ell(r)}
\right ]
=x_n+q+q^2+\cdots+q^{2n-1}t_n.
\]
Recall that in the case of $B_n,~R_{n-1}=R_{n-1}^0\cup R_{n-1}^1$, ~~where\\
$R_{n-1}^0:=\{1, s_{n-1}, s_{n-1} s_{j-2}, \dots, s_{n-1}\cdots s_1\}$ ~~and\\
$R_{n-1}^1:=\{s_{n-1}\cdots s_0,\; s_{n-1}\cdots s_0 s_1,\; 
s_{n-1}\cdots s_0 s_1 s_2,\;
 \dots,\;
s_{n-1}\cdots s_0 s_1 \cdots s_{n-1}\}$.\\
The only $r=w_{n-1}$ in $R_{n-1}$ with $K_0(r)\ne\phi$ is $r=1$, hence
the contribution of $x_n$. Similarly, the 
only $r=w_{n-1}$ in $R_{n-1}$ with $K_1(r)\ne\phi$ is
$r=s_{n-1}\cdots s_0 \cdots  s_{n-1}$ -- of length $2n-1$ -- hence the
contribution of $q^{2n-1} t_n$. This also explains the other summands
$q,\;q^2$, etc.

This implies the proof.

\qed

\bigskip

\begin{df}\label{df-fm}
For $\sg\in C_2\wr S_n=B_n$ it is natural to consider
the {\it sequence descent set} :
$$
\D_A (\sg):=\{ 0\le i\le n-1|\ \sg (i)>\sg(i+1)\}
$$
and the {\it sequence major index}
$$
\maj_A(\sg):=\sum\limits_{i\in \D_A(\sg)} i .
$$
Let 
$$
neg(\sg):=\#\{i| \sg(i)<0\}
$$
and 
define the {\it flag major index} as
$$
\fm(\sg):= 2\cdot \maj_A(\sg) + neg(\sg).
$$
\end{df}

The flag major index was introduced in~\cite{AR} 
in order to extend MacMahon classical equi-distribution theorem to $B_n$.
For a unified definition of the classical major index and the flag-major index as a length of a distinguished canonical expression see~\cite[Theorem 3.1]{AR}. 
The flag-major index has many other combinatorial and algebraic properties
which are shared with the classical major index on $S_n$,
see, for example, \cite{ABR2, AGR, HLR} and references therein.

The following theorem is a flag-major index analogue of Theorem~\ref{equid00}.

\begin{thm}\label{equid01}
For every positive integer $n$
$$
\sum\limits_{\sigma\in B_n}
\prod\limits_{i\in \M_0(\sg)}x_i 
\cdot \prod\limits_{i\in \M_{1}(\sg)}t_i
\cdot q^{\fm(\sigma)}
= 
$$
$$
(x_1+qt_1)(x_2+q+q^2+q^3t_2)\cdots (x_n+q+q^2+\cdots+q^{2n-1}t_n).
$$
\end{thm}

To prove this theorem we need the following lemma.

\begin{lem}\label{shuff-B_n}
For every $\bar \sg\in B_{n-1}$
$$
\sum\limits_{r\in R_{n-1}} q^{\fm(\bar\sg r)}=
q^{\fm(\bar\sg)}\cdot (1+q+\cdots+q^{2n-1}).
$$
\end{lem}

\noindent{\bf Proof.}
By the definition of $\fm$ (Definition~\ref{df-fm}),
$$
\sum\limits_{r\in R_{n-1}} q^{\fm(\bar\sg r)}=
\sum\limits_{r\in R_{n-1}} q^{2 maj_A(\bar\sg r)+neg(\bar\sg r)}=
$$
$$
=\sum\limits_{r\in R^0_{n-1}} q^{2 maj_A(\bar\sg r)+neg(\bar\sg r)}+
\sum\limits_{r\in R^1_{n-1}} q^{2 maj_A(\bar\sg r)+neg(\bar\sg r)}=
$$
$$
\sum\limits_{r\in R^0_{n-1}} q^{2 maj_A(\bar\sg r)+neg(\bar\sg)}
+
\sum\limits_{r\in R^1_{n-1}} q^{2 maj_A(\bar\sg r)+neg(\bar\sg)+1}=
$$
$$
=q^{neg(\bar\sg)}\cdot
\left [\sum\limits_{r\in R^0_{n-1}} q^{2 maj_A(\bar\sg r)}
+q \sum\limits_{r\in R^1_{n-1}} q^{2 maj_A(\bar\sg r)}\right ]
$$
By a theorem of Garsia and Gessel~\cite[Theorem 3.1]{GG},
$$
\sum\limits_{r\in R^0_{n-1}} q^{2 maj_A(\bar\sg r)}=
\sum\limits_{r\in R^1_{n-1}} q^{2 maj_A(\bar\sg r)}=
q^{2 maj_A(\bar\sg)}\cdot (1+q^2+\cdots+q^{2(n-1)})
$$
completing the proof of the lemma.

\qed 

\medskip

\noindent{\bf Proof of Theorem~\ref{equid01}.}
Again, by induction on $n$.
Obviously, theorem holds for $n=1$. 

Recall the definition of $r_{u,n-1}$ from Proposition~\ref{uu2}.
Then for every $\bar\sg\in B_{n-1}$, 
\begin{equation}\label{7.0}
\fm(\bar\sg r_{1,n-1} )=\fm(\bar\sg)+2n-1 \qquad \fm (\bar\sg \cdot r_{0,n-1})=\fm(\bar\sg). 
\end{equation}
Combining ~(\ref{7.0}) 
with Lemma~\ref{shuff-B_n} implies
\begin{equation}\label{7.1}
\sum\limits_{r\in R_{n-1}\setminus\{r_{0,n-1},r_{1,n-1}\}} q^{\fm(\bar\sg r)}= q^{\fm(\bar\sg)}\cdot (q+\cdots+q^{2n-2}). 
\end{equation}
Clearly, the l.h.s. in the theorem equals 
\[
\sum\limits_{\bar\sigma\in B_{n-1}} 
\sum_{r\in R_{n-1}\setminus\{r_{0,n-1},r_{1,n-1}\}} 
\prod\limits_{i\in \M_0(\bar\sigma r)}x_i 
\cdot \prod\limits_{i\in \M_{1}(\bar\sigma r)}t_i 
\cdot q^{\fm(\bar\sigma r)}
+
\]
\[
\sum\limits_{\bar\sigma\in B_{n-1}} 
\sum_{r\in \{r_{0,n-1},r_{1,n-1}\}}
\prod\limits_{i\in \M_0(\bar\sigma r)}x_i 
\cdot \prod\limits_{i\in \M_{1}(\bar\sigma r)}t_i =
\cdot q^{\fm(\bar\sigma r)}
\]
By Corollary~\ref{uu4} and (\ref{7.1}), the first sum equals
\[
\sum\limits_{\bar\sigma\in B_{n-1}} 
\sum_{r\in R_{n-1}\setminus\{r_{0,n-1},r_{1,n-1}\}} q^{\fm(\bar\sigma r)}
\cdot \prod\limits_{i\in \M_0(\bar\sigma)}x_i 
\cdot \prod\limits_{i\in \M_{1}(\bar\sigma)}t_i =
\]
\[
\cdot \prod\limits_{i\in \M_0(\bar\sigma)}x_i 
\cdot \prod\limits_{i\in \M_{1}(\bar\sigma)}t_i 
\cdot q^{\fm(\bar\sigma)}(q+\cdots+q^{n-2}).
\]
By Corollary~\ref{uu4} and (\ref{7.0}), the second sum equals
\[
\cdot \prod\limits_{i\in \M_0(\bar\sigma)}x_i 
\cdot \prod\limits_{i\in \M_{1}(\bar\sigma)}t_i 
\cdot q^{\fm(\bar\sigma)}(x_n+t_n q^{2n-1})
\]
completing the proof.

\qed


\medskip

We deduce

\begin{cor}\label{equid0}
For every positive integer $n$
$$
\sum\limits_{\sigma\in B_n}
\prod\limits_{i\in \M_0(\sg)}x_i 
\cdot \prod\limits_{i\in \M_{1}(\sg)}t_i 
\cdot q^{\ell(\sigma)}
=
\sum\limits_{\sigma\in B_n}
\prod\limits_{i\in \M_0(\sg)}x_i 
\cdot \prod\limits_{i\in \M_{1}(\sg)}t_i 
\cdot q^{\fm(\sigma)}. 
$$
%
%
Equivalently, 
for every positive integer $n$ and every pair of disjoint subsets
$B_1,B_2\subseteq \{1,\dots,n\}$
$$
\sum\limits_{\{\sigma\in B_n|\ \M_{1}(\sg)=B_1, \M_0(\sg)=B_2\}}
q^{\ell(\sigma)} =
\sum\limits_{\{\sigma\in B_n|\ \M_{1}(\sg)=B_1, \M_0(\sg)=B_2\}}
q^{\fm(\sigma)} .
$$
\end{cor}

\noindent{\bf Proof.} Combine Theorem~\ref{equid00} with
Theorem~\ref{equid01}.

\qed

\begin{cor}\label{equid2}
For every positive integer $n$
$$
\sum\limits_{\sigma\in B_n}q^{\ell(\sigma)} x^{\der_0(\sg)} t^{\der_1(\sg)}=
\sum\limits_{\sigma\in B_n}q^{\fm(\sigma)} x^{\der_0(\sg)} t^{\der_1(\sg)}=
$$
$$
(x+qt)(x+q+q^2+q^3t)\cdots (x+q+q^2+\cdots+q^{2n-1}t).
$$
\end{cor}

\noindent{\bf Proof.}
Substitute $x_1=\cdots=x_n=x$ and $t_1=\cdots=t_n=t$ in the r.h.s.
of Theorems~\ref{equid00} and~\ref{equid01}.

\qed


\section{Generalized binomial-Stirling Numbers}\label{s0}

In this section we present the generalized 
binomial-Stirling numbers,
defined by a natural recurrence relation. 

\subsection{The Recurrence: Main Examples}

\begin{df}\label{PS}
\rm{Fix three integers $a,d,r\in Z$ and let
$h(n,k)=h_{a,d,r}(n,k)$ be the numbers determined
by the following recurrence:}
\begin{eqnarray}\label{ex1}
h(n,k)=(an+dk-r)\cdot h(n-1,k)+h(n-1,k-1),
\end{eqnarray}
\rm{where $h(0,0)=1$
and $h(n,k)=0$ if $k<0$ or $n<k$.
\noindent
We call $h_{a,d,r}(n,k)$ {\it the $(a,d,r)$--binomial--Stirling numbers}.}
\end{df}

The following examples justify that terminology.

\begin{exa}\label{e1}
\rm{The three main examples of such system of numbers are the binomial 
coefficients and the two types of the Stirling numbers.}
\begin{enumerate}
\item
$a=d=0,\quad r=-1$, 
so $h(n,k)=h(n-1,k)+h(n-1,k-1)$. In this case
$h(n,k)={n\choose k}$ are the binomial coefficients.
\item
$a=r=1,\quad d=0$, so $h(n,k)=(n-1)\cdot h(n-1,k)+h(n-1,k-1)$.
Thus $h(n,k)=c(n,k)$ are the signless Stirling numbers of the first kind.
\item
$a=r=0,\quad d=1$, hence $h(n,k)=k\cdot h(n-1,k)+h(n-1,k-1)$.
Here $h(n,k)=S(n,k)$ are the Stirling numbers of the second kind.
\end{enumerate}
\end{exa}

\subsection{Matrix Product Decomposition}\label{s1}

We need to introduce some notations.
Denote the $(i,j)$-th binomial coefficients by
\[
b(i,j):=\left (\begin{array}{c}
i\\j
\end{array} \right )
\]

{\bf Notation}. We follow~\cite{ECI}.
For $1\le k\le n$,
the signless Stirling numbers of the first
kind are denoted by $c(n,k)$, $s(n,k)=(-1)^{n-k}c(n,k)$ are the Stirling numbers of the first
kind, and $S(n,k)$ denote the Stirling numbers of the second kind.
\bigskip

Let $a,d,r\in Z$ and denote $r_1=r+d$.  Assume
$a,d,r_1\ne 0$.
For the cases where some of these integers are zero see 
Remark~\ref{801}, Corollary~\ref{802} 
and Appendix 2 below.

For a positive integer $n$
construct the following $n\times n$ lower--triangular matrices:
\begin{enumerate}
\item
$c_n=(c(i,j)\mid 1\le i,j\le n)$,
\item
$s_n=(s(i,j)\mid 1\le i,j\le n)$,
\item
$S_n=(S(i,j)\mid 1\le i,j\le n)$,
\item
$P_n=(b(i,j)\mid 0\le i,j\le n-1)$,
\item
$J_n=\mbox{diag}\left (1,-1,1,-1,\ldots, (-1)^{n-1}\right )$,
\item
$a_n=\mbox{diag}\left (1,a,a^2,\ldots, a^{n-1}\right )$,
\item
$d_n=\mbox{diag}\left (1,d,d^2,\ldots, d^{n-1}\right )$,
\item
$\hat r_n=\mbox {diag}\left (1,r_1^2,\ldots, r_1^{n-1}\right )$,
where $r_1=r+d$.
\end{enumerate}
The following properties are either obvious or well known.

\begin{lem}\label{pp1}
\begin{enumerate}
\item
$J_n=J_n^{-1}$.
\item
$P_n^{-1}=J_nP_nJ_n=((-1)^{i-j}b(i,j)\mid 0\le i,j\le n-1)$.
\item
$s_n=J_nc_nJ_n$ and $S_n=s_n^{-1}$, hence $c_n=J_nS_n^{-1}J_n$.
\item
$a_ns_na_n^{-1}=\left (a^{i-j}s(i,j)\mid 1\le i,j\le n\right )$,\\

\smallskip

$d_nS_nd_n^{-1}=\left (d^{i-j}S(i,j)\mid 1\le i,j\le n\right )\quad$ and\\

\smallskip

$\hat r_nP_n \hat r_n^{-1}=
\left ( r_1^{i-j}b(i,j)\mid 0\le i,j\le n-1\right )$.
\item
The matrices $a_n,d_n,\hat r_n$ and $J_n$ commute with each other.
\item
$\lim_{a\to 0}a_ns_na_n^{-1}=
\lim_{d\to 0}d_nS_nd_n^{-1}=
\lim_{r_1\to 0}\hat r_nP_n \hat r_n^{-1}=I_n$, where $I_n$ is the $n\times n$
identity matrix.
\end{enumerate}
\end{lem}
%
%
%
%
%
%
%
\begin{thm}\label{pt1}
Let $a,d,r\in Z$, $r_1=r+d$ and $a,d,r_1\ne 0$. Let $h(n,k)$ be a system
of numbers such that the matrices $h_n=(h(i,j)\mid 0\le i,j\le n-1)$ are
lower triangular -- for all $n$.

Then $h(n,k)$ satisfy the recurrence~(\ref{ex1}) 
if and only if the following matrix equations hold for all $n$:
\begin{eqnarray}\label{pp4}
h_n=(a_nc_na_n^{-1})({\hat r}_nP_n^{-1}{\hat r}_n\,^{-1})(d_nS_nd_n^{-1}).
\end{eqnarray}
\end{thm}

\noindent{\bf Proof}.
Let $\bar h(i,j)$ denote the entries on the right-hand-side
of~(\ref{pp4}).
It suffices to show that the numbers $\bar h(n,k)$
satisfy the recurrence~(\ref{ex1}). 
%
%

\medskip

Since $\bar h(p,q)$ are given by
r.h.s.(\ref{pp4}), by matrix multiplication,
\begin{eqnarray}
\bar h(p,q)=\sum_{p\ge j\ge i\ge q}a^{p-j} (-r_1)^{j-i} d^{i-q}
c(p+1,j+1)\cdot b(j,i)\cdot S(i+1,q+1)\\ \label{pp5}
=\sum_{\infty\le i,j\le \infty}a^{p-j} (-r_1)^{j-i} d^{i-q}
c(p+1,j+1)\cdot b(j,i)\cdot S(i+1,q+1).
\end{eqnarray}
The last equality follows from the defining conditions
$c(n,k)=0$ for $k<0$ and $k>n$, and similarly
for $b(n,k)$ and $S(n,k)$. Writing the sum in this form allows us
to ignore the sum limits.

Clearly, $\bar h(0,0)=1$.
Apply now Equation~(\ref{pp5}) to show that
the numbers $\bar h(n,k)$ satisfy the
recurrence~(\ref{ex1}), 
namely, that
\begin{eqnarray}\label{pp6}
\bar h(n,k)=(an+dk-r)\cdot \bar h(n-1,k)+\bar h(n-1,k-1),
\end{eqnarray}
and this will prove the theorem.

\medskip

By~(\ref{pp5}), since $r_1=r+d$, the right hand side of~(\ref{pp6}) is:
\[
(a(n-1)+d(k+1)+(a-r_1))\cdot \bar h(n-1,k)+\bar h(n-1,k-1)
=M_{1}+M_{2} +M_3+M_4+M_5
\]
where
\[
M_1=a(n-1)\cdot \bar h(n-1,k)=~~~~~~~~~~~~~~~~~~~~~~~~~~~~~~~~~~~~~~~~~~~~~~~~~~~
\]
\begin{eqnarray}\label{n1}
~~~=\sum_{i,j}a^{n-j} (-r_1)^{j-i} d^{i-k}(n-1)\cdot
c(n,j+1)\cdot b(j,i)\cdot S(i+1,k+1),
\end{eqnarray}
\[
M_2=d(k+1)\cdot \bar h(n-1,k)=~~~~~~~~~~~~~~~~~~~~~~~~~~~~~~~~~~~~~~~~~~~~~~~~~~~
\]
\begin{eqnarray}\label{n2}
~~~=\sum_{i,j}a^{n-1-j} (-r_1)^{j-i} d^{i-k+1}
c(n,j+1)\cdot b(j,i)\cdot (k+1)\cdot S(i+1,k+1),
\end{eqnarray}
\[
M_3=\bar h(n-1,k-1)=~~~~~~~~~~~~~~~~~~~~~~~~~~~~~~~~~~~~~~~~~~~~~~~~~~~~~~~~~~~~~~
\]
\begin{eqnarray}\label{n3}
~~~=\sum_{i,j}a^{n-1-j} (-r_1)^{j-i} d^{i-k+1}
c(n,j+1)\cdot b(j,i)\cdot S(i+1,k),
\end{eqnarray}
\[
M_4=a\cdot \bar h(n-1,k)=~~~~~~~~~~~~~~~~~~~~~~~~~~~~~~~~~~~~~~~~~~~~~~~~~~~~~~~~~~~~~~
\]
\begin{eqnarray}\label{n4}
~~~=\sum_{i,j}a^{n-j} (-r_1)^{j-i} d^{i-k}
c(n,j+1)\cdot b(j,i)\cdot S(i+1,k+1)
\end{eqnarray}
and
\[
M_5=(-r_1)\cdot \bar h(n-1,k)=~~~~~~~~~~~~~~~~~~~~~~~~~~~~~~~~~~~~~~~~~~~~~~~~~~~~~~~~~~~~~~
\]
\begin{eqnarray}\label{n5}
\sum_{i,j}a^{n-1-j} (-r_1)^{j+1-i} d^{i-k}
c(n,j+1)\cdot b(j,i)\cdot S(i+1,k+1).
\end{eqnarray}
The recurrence  $(k+1)S(i+1,k+1)+S(i+1,k)=S(i+2,k+1)$ implies that
\begin{eqnarray}\label{n23}
M_{2}+M_3=\sum_{i,j}a^{n-1-j} (-r_1)^{j-i} d^{i-k+1}
c(n,j+1)\cdot b(j,i)\cdot S(i+2,k+1).
\end{eqnarray}
Replacing $i+1$ by $i$ and $j+1$ by $j$, Equation~(\ref{n23}) implies that
\begin{eqnarray}\label{pp8}
M_{2}+M_3=\sum_{i,j}a^{n-j} (-r_1)^{j-i} d^{i-k}
c(n,j)\cdot b(j-1,i-1)\cdot S(i+1,k+1).
\end{eqnarray}
Replacing $j+1$ by $j$ in~(\ref{n5}) yields
\begin{eqnarray}\label{pp7}
M_5=\sum_{i,j}a^{n-j} (-r_1)^{j-i} d^{i-k}
c(n,j)\cdot b(j-1,i)\cdot S(i+1,k+1).
\end{eqnarray}
Since $b(j-1,i)+b(j-1,i-1)=b(j,i)$, by~(\ref{pp8}) and~(\ref{pp7})
\begin{eqnarray}\label{pp9}
M_{2}+M_3+M_5=\sum_{i,j}a^{n-j} (-r_1)^{j-i} d^{i-k}
c(n,j)\cdot b(j,i)\cdot S(i+1,k+1).
\end{eqnarray}
Clearly
\begin{eqnarray}\label{pp10}
M_{1}+M_{4}=\sum_{i,j}a^{n-j} (-r_1)^{j-i} d^{i-k}n\cdot
c(n,j+1)\cdot b(j,i)\cdot S(i+1,k+1).
\end{eqnarray}
Since $n\cdot c(n,j+1)+c(n,j)=c(n+1,j+1)$, by~(\ref{pp9}),~(\ref{pp10})
and~(\ref{pp5})
we finally get
\[
M_1+M_2+M_3+M_4+M_5=~~~~~~~~~~~~~~~~~~~~~~~~~~~~~~~~~~~~~~~~~~~~~~~~~~~~~~~~~~
\]
\[
~~~~~~=\sum_{i,j}a^{n-j} (-r_1)^{j-i} d^{i-k}
c(n+1,j+1)\cdot b(j,i)\cdot S(i+1,k+1)=\bar h (n,k).~~~~~~~~~~
\]
This completes the proof.
\qed

\begin{rem}\label{801}
Theorem~\ref{pt1} 
holds when $a=0$: in that case the factor $a_nc_na_n^{-1}$ is canceled from~(\ref{pp4}) since
$\lim_{a\to 0}a_ns_na_n^{-1}=I_n$. Similarly, it holds 
when $d=0$ and when $r_1=0$.
\end{rem}

\begin{rem}
{\rm If one reverses the order in~(\ref{pp4}), it seems that the numbers
given by}
\[
h^*_n=(d_nS_nd_n^{-1})({\hat r}_nP_n^{-1}{\hat r}_n\,^{-1})(a_nc_na_n^{-1})
\]
{\rm satisfy no (obvious) recurrence.}
\end{rem}

\subsection{The Dual System $Q_{a,d,r,}(n,k)$}\label{ssec}

By a trivial induction on $i$, $h_{a,d,r}(i,i)=1$ for all $i$.
Also, by definition, the matrix
$ ((-1)^{i-j}h_{a,d,r}(i,j))_{0\le i,j\le n-1}$
is lower triangular, hence invertible. By
inversion we obtain the {\it dual system} $Q_{a,d,r}(n,k)$ :
\[
(Q_{a,d,r}(i,j))_{0\le i,j\le n-1}:=((-1)^{i-j}h_{a,d,r}(i,j))^{-1}_{0\le i,j\le n-1}
\]
and $Q_{a,d,r}(n,k)$ might be called
the $(a,d,r)$--binomial--Stirling numbers of the second kind.
%
%

\begin{df}\label{defi1}
Again, let $a,d,r\in Z$ and define $h(n,k)=h_{a,d,r}(n,k)$
via either~(\ref{ex1}) 
or~(\ref{pp4}).
\begin{enumerate}
\item
Call $h(n,k)$ {\it the
$(a,d,r)$--signless binomial--Stirling numbers of the first kind}.
Let
\[
q(n,k)=(-1)^{n-k}h(n,k),
\]
and call $q(n,k)$ {\it the $(a,d,r)$--binomial--Stirling
numbers of the first kind}.
Finally, denote
\[
h_n=(h(i,j)\mid 0\le i,j\le n-1 ),\quad\mbox{and}\quad
q_n=(q(i,j)\mid 0\le i,j\le n-1 ),
\]
$n\times n$ lower--triangular matrices.
\item
Define the numbers $Q(i,j)$ by inverting the matrices $q_n$:
\[
Q_n=(Q(i,j)\mid 0\le i,j\le n-1 )=q_n^{-1}=
(q(i,j)\mid 0\le i,j\le n-1 )^{-1}=~~~~~~~~~~~
\]
\[
~~~~~~~~~~~~~~~~~~~~~~~~~~~~~~~~=((-1)^{i-j}h(i,j)\mid
0\le i,j\le n-1 )^{-1}.
\]
Also $Q(i,j)=0$ if $j<0$ or if $i<j$.
The definition of $Q(i,j)$ is independent of $n$, provided $i\le n$.

Call $Q(n,k)=Q_{a,d,r}(n,k)$ {\it the $(a,d,r)$--binomial--Stirling
numbers of the second kind}.
\end{enumerate}
\end{df}

The following theorem shows that such binomial--Stirling
numbers of the second kind are just a binomial--Stirling system, but with
$(d,a,r+d-a)$ replacing $(a,d,r)$.

\medskip

Clearly, $q_n=J_n h_n J_n$, hence $ h_n=J_n Q_n^{-1} J_n$.
With the notations of Subsection~\ref{s1} we have
\begin{thm}\label{c1}
Let $Q(n,k)=Q_{a,d,r}(n,k)$ denote the
$(a,d,r)$-Stirling numbers of the second kind,
with corresponding matrices $Q_n$, then
\begin{enumerate}
\item
\begin{eqnarray}\label{s2}
Q_n=(d_nc_nd_n^{-1})({\hat r}_nP_n^{-1}{\hat r}_n\,^{-1})(a_nS_na_n^{-1}).
\end{eqnarray}
\item
The numbers $Q(n,k)=Q_{a,d,r}(n,k)$ satisfy the following recurrence, which is
 ``dual'' to the recurrence~(\ref{ex1}) :
\begin{eqnarray}\label{s3}
Q(n,k)=(dn+ak-r_2)\cdot Q(n-1,k)+Q(n-1,k-1),
\end{eqnarray}
where $r_2=r+d-a$. Thus $Q(n,k)=Q_{a,d,r}(n,k)=h_{d,a,r+d-a}(n,k)$.
\end{enumerate}
\end{thm}
\noindent{\bf Proof.}
1.  By  Lemma~\ref{pp1}.1 and Definition~\ref{defi1},
$h_n^{-1}=J_n Q_nJ_n$.
Inverting~(\ref{pp4}) implies that
\[
J_n Q_n J_n=
(d_nS_n^{-1}d_n^{-1})({\hat r}_nP_n{\hat r}_n\,^{-1})(a_nc_n^{-1}a_n^{-1}).
\]
Applying Lemma~\ref{pp1}, deduce that
\[
Q_n =
J_n(d_nJ_nc_nJ_nd_n^{-1})({\hat r}_nP_n{\hat r}_n\,^{-1})
(a_nJ_nS_nJ_na_n^{-1})J_n=~~~~~~~~~~~~~~~~~
\]
\[
~~~~~~~~~~~~~~~=(d_nc_nd_n^{-1})({\hat r}_nJ_nP_nJ_n{\hat r}_n\,^{-1})
(a_nS_na_n^{-1}),
\]
and the proof follows by Lemma~\ref{pp1}.2.

\medskip

2. By Theorem~\ref{pt1},~(\ref{s3}) and~(\ref{s2}) are equivalent,
with $r_1=r_2+a$.

\qed
%
%

\begin{cor}\label{802} 
In particular the two systems, of $h_{a,0,1}(n,k)$ and of 
$h_{0,a,1-a}(n,k)$, are dual to each other: the systems of $h_{0,a,1-a}(n,k)$
is obtained by inverting the corresponding lower-triangular
matrix with entries $(-1)^{n-k}h_{a,0,1}(n,k)$.
\end{cor}

\begin{rem}\label{lamt1}
Theorem~\ref{c1} holds when either of $a$,  $d$ or 
$r_1=r+d$ equal zero.
%
%
For example, when $a=0$ apply $\lim_{a\to 0}$ to~(\ref{s2}), using 
Lemma~\ref{pp1}.6, see Remark~\ref{801}. 
Similarly, for $d=0$ or $r_1=0$.
\end{rem}

\medskip

\subsection{The $(a,d,r,\ell)$ Systems}
Let $\ell=\ell'\ne 0$ and consider the system of numbers 
$g(n,k)=g_{a',d',r',\ell'}(n,k)$ given by the following
$(a',d',r',\ell')$--recurrence:

\noindent Again $g(0,0)=1$, and 

\begin{eqnarray}\label{s33}
g(n,k)=(a'n+d'k-r')\cdot g(n-1,k)+\ell'\cdot g(n-1,k-1).
\end{eqnarray}
By a trivial induction one proves:

\begin{rem}\label{lo1}

 Let $a=a'/\ell,~d=d'/\ell,~r=r'/\ell$, $\ell'=\ell$, and let
$h(n,k)=h_{a,d,r}(n,k)$ be given as in Equation~(\ref{ex1}). 
Then for all $n$ and $k$
\begin{eqnarray}\label{s333}
g_{a',d',r',\ell'}(n,k)=\ell^n\cdot h_{a,d,r}(n,k).
\end{eqnarray}
Thus 
\begin{eqnarray}\label{sch3}
g_{a,d,r,1}(n,k)=h_{a,d,r}(n,k).
\end{eqnarray}
\end{rem}

\medskip

Similar to the dual system $Q(n,k)=Q_{a,d,r}(n,k)$ of $h_{a,d,r}(n,k)$,
construct the dual
system $V(n,k)=V_{a',d',r',\ell'}(n,k)$ of $g_{a',d',r',\ell'}(n,k)$ 
as follows:

Let $v(n,k)=(-1)^{n-k}g(n,k)$,~$v_n=[v(i,j)\mid 1\le i,j\le n]$
and the numbers $V(n,k)=V_{a',d',r',\ell'}(n,k)$ are given by the matrix
equation
$v_n^{-1}=[V(i,j)\mid 1\le i,j\le n]$.

\medskip

By matrix inversion one proves 
\begin{rem}\label{lo2}
For all $n$ and $k$
\begin{eqnarray}\label{s3333}
V_{a',d',r',\ell'}(n,k)=\ell^{-k}Q_{a,d,r}(n,k).
\end{eqnarray}
\end{rem}

\section{Realizations of the dual systems}\label{col9}

%
%

In Sections~\ref{col5} and~\ref{col6} two systems of 
binomial-Stirling numbers are realized by certain statistics on 
colored permutations. It is shown here that
these two systems are dual to each other - in 
the sense of Section~\ref{s0}.

\begin{rem}\label{lo3}

\noindent 1.
Note that Corollary~\ref{ch2} can be considered as 
a ``wreath-product-realization'' of the 
system $g_{a,0,\ell,\ell}(n,k)$ with $0\le \ell\le a-1$:
the recurrence of $g_L(n,k) $ there implies that 
$g_L(n,k)=g_{a,0,\ell,\ell}(n,k)$, thus
\[
g_{a,0,\ell,\ell}(n,k)=\#\{\sg\in C_a\wr S_n\mid\; \m _L (\sg)=k\}.
\]
In particular, if 
$L=\{u\}$ then $\ell =1$, $g_{a,0,1 ,1}(n,k)=h_{a,0,1}(n,k)$
and we have
\[
h_{a,0,1}(n,k)=\#\{\sg\in C_a\wr S_n\mid\; \m _u (\sg)=k\}.
\]

\noindent 2.
Similarly, Theorem~\ref{sec8} (with $d$ replacing $a$)
is a ``wreath-product-realization'' of the 
system $g_{0,d,\ell -d ,\ell}(n,k)$ with $0\le \ell\le d-1$:
the recurrence of $f_L(n,k) $ there implies that 
$f_L(n,k)=g_{0,d,\ell -d ,\ell}(n,k) $, and by Definition~\ref{ws2}
\[
g_{0,d,\ell -d ,\ell}(n,k)=\#\{\sg\in  C_d\wr S_n\mid des_L(\sg)=\m_L(\sg)=k\}.
\]
In particular, if 
$L=\{u\}$ then $\ell =1$, $g_{0,d,\ell -d ,\ell}(n,k)=h_{0,d,1-d}(n,k)$, 
hence
\[
h_{0,d,1-d}(n,k)=\#\{\sg\in  C_d\wr S_n\mid des_u(\sg)=\m_u(\sg)=k\}.
\]
\noindent 3. 
Summing the above on $k$ implies
\[
\sum_k
g_{0,d,\ell -d ,\ell}(n,k)=\#\{\sg\in  C_d\wr S_n\mid des_L(\sg)=\m_L(\sg)\}.
\]
\end{rem}

\medskip

This leads to the $(d,r)-Bell$ numbers and with 
the following wreath-product-realization:

\begin{df}\label{bn1}
Recall the numbers $h_{0,d,r}(n,k)=g_{0,d,r,1}(n,k)$, denote
\[
b_{d,r}(n)=\sum_k h_{0,d,r}(n,k)=\sum_k g_{0,d,r,1}(n,k)
\]
and call these the (d,r)-Bell numbers.
\end{df}
Note that by Example~\ref{e1}.3 $h_{0,1,0}(n,k)=S(n,k)$ are the
Stirling numbers of the second kind, therefore $b_{1,0}(n)$ are the (ordinary)
Bell--numbers. Further properties of the 
$(d,r)$--Bell numbers are given in Appendix 2.

By Remark~\ref{lo3}.3 and the above definition,

\begin{cor}\label{bn2}
\[
b_{d,1-d}(n)=
\# \{\sigma\in C_d\wr S_n\ |\ des_u(\sigma)=\ \der_u (\sigma)\}.
\]
\end{cor}

Recall from~\cite[Propositions 10.8 and 10.10]{RR2} 
that the signed Stirling number of the first kind, $s(n,k)$, is equal to
$$
(-1)^{n-k}\cdot \# \{\pi\in S_n\ |\ \der (\pi)=k\},
$$
while the Stirling number of the second kind, $S(n,k)$, is equal to
$$
\# \{\pi\in S_n\ |\ des (\pi)=\ \der (\pi)=k\}.
$$
These numbers form inverse matrices,
see, e.g., \cite[Prop. 1.4.1.a]{ECI}.
This phenomenon is generalized to wreath products.
\medskip


\begin{thm}\label{dual}
For every positive integers $a$, $N$, and every subset 
$L\subseteq \{0,\dots,a-1\}$ of size $\ell$ let $s_{L,N}$ 
be the $N\times N$ matrix
whose entries are given by
$$
s_{L,N}(n,k):={(-1)^{n-k}\over \ell^n}\cdot \# \{\sigma\in C_a\wr S_n\ |\ 
\der_L(\sigma)=k\} \qquad (0\le k,n\le N)
$$
and $S_{L,N}$ be the $N\times N$ matrix whose entries are
$$
S_{L,N}(n,k):=
{1\over \ell^n}\cdot 
 \# \{\sigma\in C_a\wr S_n\ |\ des_u(\sg)=\ \der_u(\sigma)=k\}
\qquad (0\le k,n\le N).
$$
Then 
$$
S_{L,N}^{-1}= s_{L,N}.
$$
\end{thm}

%

\noindent{\bf Proof.}
First note that the results in Section~\ref{s0} hold for 
any rational (essentially real) $a,d,r$.
Thus, by Remarks~\ref{lo3}.(1) and~\ref{lo1},
$$
\ell^{-n}\cdot \#\{\sg\in C_a\wr S_n\mid \m _L (\sg)=k\}=
\ell^{-n}g_{a,0,\ell,\ell}(n,k)=
h_{a/\ell,0,1}(n,k).
$$
Similarly by Remarks~\ref{lo3}.(2) and~\ref{lo1},
$$
\ell^{-n}\cdot \# 
\{\sg\in  C_a\wr S_n\mid des_L(\sg)=\m_L(\sg)=k\}=
$$
$$
=\ell^{-n}g_{0,a,\ell-a,\ell}(n,k)=
h_{0,a/\ell,1-a/\ell}(n,k).
$$
Corollary~\ref{802} completes the proof.

\qed

\section{Appendix 1: Proof of Proposition~\ref{ddes-L}}\label{app1}

For every element
$\sg\in C_a\wr S_n$ and $L\subseteq \{0,\dots,a-1\}$ define
$$
\inv_L(\sg):=\{i<j|\ \sg(i)>_L \sg(j)\}.
$$
For $1\le i\le n$ denote  the color of $\sg(i)$ by
$z_i(\sg)$. Namely, $z_i(\sg)=j$ if 
$\sg(i)=\al^j |\sg(i)|$.

We will apply the following combinatorial formula for the length function.

\begin{lem}\label{Bagno2} \cite[Theorem 3.3.3]{B}
For every positive integers $a$ and $n$, and every element
$\sg\in C_a\wr S_n$
$$
\ell(\sg)= \inv_{\bar 0}(\sg)+\sum\limits_{\sg(i)<_{\bar 0} 0} (|\sg(i)|-1) + \sum\limits_{j=1}^n z_j(\sg),
$$
where $\bar 0:=\{1,\dots,a-1\}$.
\end{lem}


\begin{cor}\label{des10}
For every element
$\sg\in C_a\wr S_n$ and $0\le i\le n-1$,
$$
\ell(\sg s_i)<\ell(\sg)\Longleftrightarrow \sg(i)>_{\bar 0} \sg(i+1)
$$
where we assume $\sg(0)=0$ and $\bar 0:=\{1,\dots,a-1\}$.
\end{cor}

\noindent{\bf Proof.}  
By the definition of the order $<_{\bar 0}$ together with
Fact~\ref{b1}(1) the corollary holds for $i=0$.
For $i\ne 0$ the corollary follows from Lemma~\ref{Bagno2} 
together with Fact~\ref{b1}(2).

\qed

\noindent{\bf Proof of Proposition~\ref{ddes-L}.}
By Corollary~\ref{des10},
$$
\{ 0\le i\le n-1|\ \ell(v_u^{-1}\sg s_i)< \ell(v_u^{-1}\sg) \}= 
\{ 0\le i\le n-1|\ v_u^{-1}\sg (i)>_{\bar 0} v_u^{-1}\sg(i+1) \}. 
$$
One may easily verify that 
$$
v_u^{-1}\sg (i)>_{\bar 0} v_u^{-1}\sg(i+1)\Longleftrightarrow
\sg (i)>_{\bar u} \sg(i+1),
$$
where $\bar u:=\{0,\dots,a-1\}\setminus u$.

This proves that
for every $0\le u\le a-1$ 
\begin{equation}\label{**}
\D_{\bar u}(\sg) = 
\{
0\le i\le n-1|\ 
\ell(v_u^{-1}\sg s_i)<
\ell(v_u^{-1}\sg)
\}.
\end{equation}

\noindent
Proposition~\ref{ddes-L} is deduced from~(\ref{**})
by observing that $>_{\bar u}$
is the reverse order of $>_u$; hence $\D_u(\sg)=\{0,\dots,n-1\}\setminus \D_{\bar u}$. 


\qed

\section{Appendix 2: Further Properties of the Generalized binomial-Stirling and Bell Numbers}\label{app2}

In this appendix we study some further properties of the generalized binomial-Stirling and Bell numbers, introduced in Section~\ref{s0}.


%
\begin{pro}\label{sp3}
Let $d=0$; namely, the numbers $h_{a,0,r}(n,k)$ satisfy
the recurrence~(\ref{ex1}) 
with $d=0$.
Let
\begin{eqnarray}\label{sp1}
f_n(x):=\sum_k h_{a,0,r}(n,k)x^k.
\end{eqnarray}
Then
\begin{equation}
f_n(x)=(x+a-r)(x+2a-r)\cdots  (x+na-r).
\end{equation}
In particular
$\quad\sum_k h_{a,0,r}(n,k)=(a-r+1)(2a-r+1)\cdots (na-r+1)$.
\end{pro}
\noindent{\bf Proof.} For $n\ge 1$ let
 $\bar h(n,k)$ be the coefficient of $x^k$
in the following expansion:
\[
\bar f_n(x)=(a-r+x)(2a-r+x)\cdots (na-r+x)=\sum_k \bar h(n,k)x^k,
\]
and define $\bar h(0,0):=1$.

Then $\bar f_n(x)=(x+na-r)\bar f_{n-1}(x)=
(na-r)\bar f_{n-1}(x)+x\bar f_{n-1}(x)$. It easily follows that
$\bar h(n,k)$ satisfies the same recurrence as $h_{a,0,r}(n,k)$,
which implies that $\bar h(n,k)=h_{a,0,r}(n,k)$.

\qed

\bigskip
%

%

%
%


%
In the rest of this section we study the binomial-Stirling numbers 
with $a=0$, namely
$h(n,k)=h_{0,d,r}(n,k)$,
and deduce further results about these numbers and 
their sums, 
the $(d,r)$--Bell numbers. We follow closely Section 1.6 of~\cite{W}.

\smallskip

Denote
\begin{eqnarray}\label{sp2}
g_k(y):=\sum_nh(n,k)y^n=\sum_n h_{0,d,r}(n,k)y^n.
\end{eqnarray}
Thus $h_{0,d,r}(n,k)$ is the coefficient of $y^n$ in $g_k(y)$.


\begin{pro}\label{sp4}
Let $a=0$, namely the numbers $h(n,k)=h_{0,d,r}(n,k)$ satisfy
recurrence ~(\ref{ex1}) 
with $a=0$:
\begin{eqnarray}\label{sp5}
h(n,k)=(dk-r)\cdot h(n-1,k)+h(n-1,k-1).
\end{eqnarray}
Then
\begin{eqnarray}\label{sp6}
g_k(y)=\frac{y^{k}}{(1-(-r)y)(1-(d-r)y)\cdots (1-(kd-r)y)}.
\end{eqnarray}
\end{pro}
\noindent{\bf Proof.} Define $\bar g_k(y)$ and $\bar h(n,k)$
 via the expansion of the following ratio:
\begin{eqnarray}\label{sp7}
\bar g_k(y)=
\frac{y^{k}}{(1-(-r)y)(1-(d-r)y)\cdots (1-(kd-r)y)}
=\sum_n\bar h(n,k)y^n.
\end{eqnarray}
Clearly
\[
\bar g_k(y)=\frac{y}{1-(dk-r)y}\cdot \bar g_{k-1}(y),
\]
hence
$\bar g_k(y)=(dk-r)\cdot y\cdot \bar g_k(y)+ y\cdot\bar g_{k-1}(y)$,
namely
\[
\sum_{n} \bar h(n,k)y^n=
\sum_{n} (dk-r)\bar h(n-1,k)y^n+ \sum_{n} \bar h(n-1,k-1)y^n.
\]
Comparing coefficients, it follows that $\bar h(n,k)$ satisfy the same
recurrence~(\ref{sp5})
 as $h(n,k)$, hence $h(n,k)=\bar h(n,k)$, which completes the proof.
\qed

\begin{cor} (See~\cite[Ex. 16 in Ch. 1]{ECI})
\begin{eqnarray}\label{for3}
h_{0,d,r}(n,k)=\sum(-r)^{a_0-1}\cdot (d-r)^{a_1-1}\cdots (kd-r)^{a_{k}-1},
\end{eqnarray}
the sum being over all compositions $a_1+\cdots +a_{k+1}=n+1$ where
all $a_i\ge 1$.
\end{cor}
It should be interesting to give Equation~(\ref{for3}) a purely combinatorial
proof.

\smallskip

The following proposition extends~\cite[(1.6.7)]{W}.

\begin{pro}\label{for1}
\begin{eqnarray}\label{ffor2}
h_{0,d,r}(n,k)=\sum_{t=0}^{k}(-1)^{k-t}\cdot\frac{(td-r)^{n}}{d^k\cdot t!(k-t)!}.
\end{eqnarray}
\end{pro}
\noindent{\bf Proof.}
Let 
\[
g_k^*(y)=y^{-k}g_k(y)=\frac{1}{(1-(-r)y)(1-(d-r)y)\cdots (1-(kd-r)y)},
\]
and notice that $h_{0,d,r}(n,k)$ is the coefficient of $y^{n-k}$ in $g_k^*(y)$.
Applying partial fractions, this can be written as
\[
\frac{1}{(1-(-r)y)(1-(d-r)y)\cdots (1-(kd-r)y)}
=\sum_{t=0}^{k}\frac{\alpha_t}{1-(td-r)y}
\]
with some $\alpha_t\in R$.

To calculate $\alpha_{t}$, multiply both sides by $1-(td-r)y$,
then substitute $y= 1/(td-r)$.
On the right we get $\alpha_t$ and on the left --
\[
\frac {1}
{(1-(-r)y)\cdots (1-((t-1)d-r)y)(1-((t+1)d-r)y)\cdots (1-(kd-r)y)}=
\]

\smallskip

\[
=\frac {1}
{\left ( 1- \frac{0\cdot d-r}{t\cdot d-r}\right )\cdots
\left ( 1- \frac{(t-1)\cdot d-r}{t\cdot d-r}\right )\cdot
\left ( 1- \frac{(t+1)\cdot d-r}{t\cdot d-r}\right )\cdot
\left ( 1- \frac{k\cdot d-r}{t\cdot d-r}\right )}=
\]

\smallskip

\[
\frac{(td-r)^k}
{td\cdot (t-1)d\cdots d\cdot (-d)\cdot (-2d)\cdots (-(k-t)d)}=
\]

\smallskip

\[
=(-1)^{k-t}\frac{(td-r)^k}{d^k\cdot t!(k-t)!}.
\]
Deduce that
\[
\alpha _t=(-1)^{k-t}\frac{(td-r)^k}{d^k\cdot t!(k-t)!}.
\]

Recall that $h(n,k)$ is
the coefficient of $y^{n-k}$ in $g_k^*(y)$, and that
\[
g_k^*(y)=
\sum_{t=0}^{k}\frac{\alpha_t}{1-((td-r)y}.
\]
Thus
\[
h(n,k)=\left [ y^{n-k}\right ]
\sum_{t=0}^{k}\frac{\alpha_t}{1-(td-r)y}=
\]
\[
\sum_{t=0}^{k}\left [ y^{n-k}\right ] \frac{\alpha_t}{1-(td-r)y}=
\sum_{t=0}^{k}(td-r)^{n-k}\alpha_t=
\]
\[
=\sum_{t=0}^{k}(td-r)^{n-k}(-1)^{k-t}\cdot\frac{(td-r)^k}
{d^k\cdot t!(k-t)!}=
\]
\[
=\sum_{t=0}^{k}(-1)^{k-t}\cdot\frac{(td-r)^{n}}{d^k\cdot t!(k-t)!}.
\]
\qed
%


\medskip

Recall the $(d,r)$--Bell numbers
$b_{d,r}(n)=\sum_kh_{0,d,r}(n,k)$ from Definition~\ref{bn1}.
We have the following formula for these numbers, extending
a remarkable result of Dobinski~\cite{D}.
\begin{pro}\label{for2}
\[
b_{d,r}(n)=\frac{1}{e^{1/d}}\sum_{t=0}^{\infty}
\frac{(td-r)^n}{t!d^t}.
\]
\end{pro}

\noindent{\bf Proof 
.} We continue to follow~\cite{W}.

Choose $M$ large enough, then, by the previous proposition,
\[
b_{d,r}(n)=\sum_{k=0}^M\sum_{t=0}^k
(-1)^{k-t}\cdot\frac{(td-r)^{n}}{d^k\cdot t!(k-t)!}=
\]
\[
=\sum_{t=0}^M\frac{(td-r)^{n}}{t!d^t}\cdot\sum_{k=t}^M\frac{(-1)^{k-t}}
{(k-t)!d^{k-t}}=
\]
\[
=\sum_{t=0}^M\frac{(td-r)^{n}}{t!d^t}\cdot\sum_{s=0}^{M-t}
\frac{(-1)^{s}}{(s)!}\cdot \left (\frac{1}{d}\right )^s.
\]
The proof now follows by sending $M$ to infinity, since then, the second
factor becomes
\[
\sum_{s=0}^{\infty}
\frac{(-1)^{s}}{(s)!}\cdot \left (\frac{1}{d}\right )^s=
\frac{1}  {e^{1/d}}.
\]

\qed
%

%
%

\begin{cor}\label{bn3}
For every positive $n$
\begin{equation}\label{bn33}
\# \{\sigma\in B_n\ |\ des(\sigma)=\ \der_1 (\sigma)\}
=\frac{1}{\sqrt e}\sum_{t=0}^{\infty}
\frac{(2t+1)^n}{t!2^t},
\end{equation}
where $des(\sg)=\#\{0\le i\le n-1|\ \ell(\sg s_i)<\ell(\sg)\}$
is the standard descent number.
\end{cor}

\noindent{\bf Proof.} Combine Corollary~\ref{bn2},
with Proposition~\ref{for2} (letting $d=2$ and $r=-1$).
The identity $des(\sg)=des_1(\sg)$ ($\forall \sg\in B_n$) (see Example~\ref{ddes-L-example}.3)
completes the proof. 

\qed

\begin{df}
Let $B_{d,r}(x)$ be the exponential generating function of the $b_{d,r}(n)$'s:
\[
B_{d,r}(x)=\sum _{n=0}^{\infty} b_{d,r}(n)\frac{x^n}{n!}.
\]
\end{df}
\begin{pro}
\[
B_{d,r}(x)=\exp \left ( \frac{e^{dx}-drx-1}{d} \right ).
\]
\end{pro}
\noindent{\bf Proof.} By definition,
$\;b_{d,r}(0)=1$; hence, by Proposition~\ref{for2},
\[
B_{d,r}(x)-1=\frac{1}{e^{1/d}}\sum _{n=1}^{\infty}\frac{x^n}{n!}
\sum_{t=0}^{\infty}
\frac{(td-r)^n}{t!d^t}=
\]
\[
\frac{1}{e^{1/d}}\sum_{t=0}^{\infty}\frac{1}{t!d^t}
\sum _{n=1}^{\infty}\frac{[(td-r)x]^n}{n!}=
\]
\[
\frac{1}{e^{1/d}}\sum_{t=0}^{\infty}\frac{1}{t!d^t}\cdot
\left ( e^{(td-r)x}-1\right )=
\]
\[
-1+\frac{e^{-rx}}{e^{1/d}}\sum_{t=0}^{\infty}\frac{1}{t!d^t}\cdot
e^{tdx}.
\]
Thus
\[
B_{d,r}(x)=e^{-(drx+1)/d}\sum_{t=0}^{\infty}\frac{1}{t!}\cdot
\left (\frac{e^{dx}}{d} \right )^t=~~~~~~~~~~~~~~~~~~~~~~~~~~
\]
\[
~~~~~~~~~~~~~~~~~~~~~~~~=e^{-(drx+1)/d} \cdot e^{(e^{dx}/d)}=
\exp \left ( \frac{e^{dx}-drx-1}{d} \right ).
\]
\qed

\bigskip

\noindent{\bf Acknowledgements.}
We thank Christian Krattenthaler for some useful references.

\end{document}